\documentclass[prl,aps]{revtex4}
\usepackage{graphicx}

\begin{document}

\title{Norm of Bethe Wave Function as a Determinant}
 \bigskip

\author{Vladimir  Korepin}
\affiliation{Steklov  Mathematical  Institute of Academy of Sciences, Sankt Petersburg, Russia}
\date{Feb 26, 1981}

\begin{abstract}
This is a historical note.  Bethe Ansatz solvable models are considered, like XXZ Heisenberg anti-ferromagnet  and
Bose gas with delta interaction. Periodic boundary conditions lead to Bethe equation. The square of the norm of Bethe wave function is equal to a determinant of  linearized system of Bethe equations (determinant of matrix of second derivatives of Yang action).   The  proof  was first  published in {\it Communications in Mathematical Physics, vol 86, page 391 in  l982}.  Also {\bf domain wall boundary conditions for 6 vertex model} were discovered in the same paper [see Appendix D].
These play an important role for algebraic combinatorics: alternating sign matrices, domino tiling and plane partition.

\end{abstract}


\maketitle

\section{Introduction}
 
 Many two dimensional models have been solved by means of the Bethe Ansatz,
see for example \cite{beth,yy,lieb,bpf,mcg} \footnote{Experts believe that each universality class in 2D contain at list one solvable model.}. Quantum inverse scattering
method (QISM) \cite{fadd,fadtakh, book} discloses   algebraic nature of these solutions.
Michael  Gaudin
studied norms of Bethe wave functions  for the quantum nonlinear Schroedinger equation \footnote{the model is also known as Bose gas with delta interaction} \cite{gaudin} and
suggested  a remarkable conjecture that the norm of the eigenfunction is equal to a Jacobian.  In \cite{origin}
these formulae are proved and a more general result
is obtained. The norms are calculated for any exactly solvable models with the
R matrix either of the X X X model or of the X X Z Heisenberg models.
This note gives an idea about an approach of the paper.

 \section{Notations}
 
 First of all let us remind the reader of some notations of the QISM. Eigenfunctions
of the Hamiltonian of the physical system are constructed by means of the
monodromy matrix of an auxiliary linear problem $T(\lambda)$. In our case $T(\lambda)$ is a 2 x 2
matrix, its matrix elements being quantum operators, which depend on the
spectral parameter $\lambda$:
 \[T(\lambda)= \left( \begin{array}{cc}
A(\lambda) & B (\lambda)  \\
C (\lambda)& D (\lambda) \end{array} \right)\] 
 Commutation relations of these entries  are given by the formula:
 
 \begin{equation} \label{rtt}
R(\lambda, \mu)\left(T(\lambda) \otimes I  \right) \left( I\otimes T(\mu) \right)  = \left( I\otimes T(\mu)  \right) \left(T(\lambda) \otimes I  \right)R(\lambda, \mu)
 \end{equation}
 Here $I$ is the unit 2x2matrix, $R(\lambda, \mu)$ is a 4x4 matrix with c-number elements.
 Another way to write  this equation  is:
  \begin{equation} \label{rttind}
R(\lambda, \mu)_{\alpha , \beta}T(\lambda)_{\alpha}  T(\mu)_{\beta}  = T(\mu)_{\beta}T(\lambda)_{\alpha}R(\lambda, \mu)_{\alpha , \beta}
 \end{equation}
 This means hat $T(\lambda)_{\alpha}$ is a 2 x 2 matrix which acts on the 2 dimensional space with
index $\alpha$  and $ T(\mu)_{\beta} $ is a matrix which acts on space with index $\beta$. We assume that  $\alpha \neq \beta$. The matrix  $R(\lambda, \mu)_{\alpha , \beta}$, acts in the tensor product of these two spaces. We shall deal with R matrices of
the following form:
   \[ R(\lambda, \mu) = \left( \begin{array}{cccc}
f(\mu , \lambda) & 0 & 0 & 0  \\
0& 1 & g (\mu , \lambda) &0 \\
0 &g (\mu , \lambda) & 1& 0 \\
0&0&0&f(\mu , \lambda) \end{array} \right)\] 
For models of the XXZ type:
\begin{equation} \label{xxzel}
f(\lambda , \mu )= \frac{\sinh(\lambda -\mu +2i\eta)}{\sinh(\lambda -\mu)}, \qquad \qquad g(\lambda , \mu )= \frac{i\sin(2\eta)}{\sinh(\lambda -\mu)}
 \end{equation}
For models of  XXX type:
\begin{equation} \label{xxxel}
f(\lambda , \mu )= \frac{(\lambda -\mu +2i\kappa)}{(\lambda -\mu)}, \qquad \qquad g(\lambda , \mu )= \frac{i\kappa}{(\lambda -\mu)}
 \end{equation}
Here $\lambda$ and $\mu$ are spectral parameters and  $\eta$ and $\kappa$ are coupling constants. Note that XXX is a limiting case of XXZ. We can also write the $R$ matrix using notations of formula (\ref{rttind}).

\begin{eqnarray}  
& R(\lambda, \mu)_{\alpha , \beta} =  \cos(\eta) \frac{\sinh(\mu-\lambda  +i\eta)}{\sinh(\mu-\lambda )} I_{\alpha} I_{\beta} +i\sin(\eta)
\frac{\cosh(\mu-\lambda  +i\eta)}{\sinh(\mu-\lambda )}
\sigma_{\alpha}^{3}\sigma_{\beta}^{3}\\
&+ i \frac{\sin(2\eta)}{\sinh(\mu-\lambda )} (\sigma^{+}_{\alpha}\sigma^{-}_{\beta}+\sigma^{-}_{\alpha}\sigma^{+}_{\beta} )
\end{eqnarray} 
Here $\sigma$ are the standard Pauli matrices and $2\sigma^{\pm}=\sigma^{1}\pm\sigma^{2}$.
Let us write down some commutation relations of (\ref{rtt}) explicitly 
\begin{eqnarray} {\label{comut}}
& [B(\lambda), B(\mu)]=0, \quad    [C(\lambda), C(\mu)]=0, \quad [A(\lambda) + D(\lambda), A(\mu) + D(\mu)]=0 \\
& A(\mu)B(\lambda)= f(\mu, \lambda) B(\lambda) A(\mu) + g(\lambda , \mu) B(\mu) A(\lambda) \\
& D(\mu)B(\lambda)= f( \lambda , \mu) B(\lambda) D(\mu) + g(\mu , \lambda ) B(\mu) D(\lambda) \\
& [C(\lambda^{c}), B(\lambda^{b})]= g(\lambda^{c} , \lambda^{b}) \{ A(\lambda^{c}) D (\lambda^{b}) - A(\lambda^{b}) D (\lambda^{c}) \}
\end{eqnarray} 
Other commutation relations can be found in the book \cite{book}.
 Pseudo-vacuum $|0>$ and dual pseudo-vacuum $<0|$ are important:
 \begin{eqnarray} {\label{vac}}
 &C(\lambda)|0>=0, \qquad A(\lambda)|0>=a(\lambda)|0>, \qquad D(\lambda)|0>=d(\lambda)|0> \\
 &<0|B(\lambda)=0, \qquad <0| A(\lambda) =<0| a(\lambda), \qquad <0| D(\lambda) =<0| d(\lambda)
 \end{eqnarray} 
 Here $a(\lambda)$ and $d(\lambda)$ are complex valued functions. The space in which the operators $A(\lambda)$, $B(\lambda) $, $C((\lambda)$ and $D(\lambda)$ act is constructed in \cite{book}.
The norms and scalar products
in question are functionals of these $a(\lambda)$ and $d(\lambda)$. We shall vary $a(\lambda)$  and $d(\lambda)$ and
study the dependence of scalar products on these functional arguments.
The Hamiltonian of the physical system in question is expressed in terms of
the transfer matrix $t(\mu)=A(\lambda)+D(\lambda)$. The Hamiltonian and the transfer matrix
have common eigenfunctions which are constructed as follows. Put
\begin{equation} \label{wf}
\psi_{N}(\{ \lambda_{j} \})=B(\lambda_{1})\ldots B(\lambda_{N}) |0>
\end{equation}
and suppose that $ \lambda_{j}$ satisfy the system of the transcendental equations (TE) :
\begin{equation} \label{be}
\frac{a(\lambda_{n})}{d(\lambda_{n})} \prod_{j=1}^{N}\frac{f(\lambda_{n} , \lambda_{j} )}{f(\lambda_{j} , \lambda_{n} ) }   =1 \qquad \mbox{here} \qquad j\neq n
\end{equation}
Then $\psi_{N}(\{ \lambda_{j} \})$ is an eigenfunction of $t(\lambda)$ with the eigenvalue:
\begin{equation} \label{ev}
\theta (\mu)= a(\mu)\prod_{j=1}^{N}f(\mu , \lambda_{j}) + d(\mu)\prod_{j=1}^{N}f( \lambda_{j}, \mu) 
\end{equation}
Here $N$ is called the number of particles.
Note that
\begin{equation} \label{dwf}
\overline{\psi_{N}(\{ \lambda_{j} \})}=<0|C(\lambda_{1})\ldots C(\lambda_{N}) 
\end{equation}
is a dual eigenfunction for $t(\mu)$ with the same eigenvalue. Pauli principle was proved in the original publication \cite{origin}, see also  \cite{book}.  We take all $\lambda_{j}$ to be different.
Finally let us present two remarks. First of all new variables 
\begin{equation} \label{lhs}
\phi_{k}= i \ln \frac{a(\lambda_{k})}{d(\lambda_{k})} + i\sum_{j=1, j\neq k}^{N} \ln \frac{f(\lambda_{k}, \lambda_{j})} {f(\lambda_{j}, \lambda_{k})}, \qquad k=1, \ldots N
\end{equation}
are convenient. For example (\ref{be}) can be written as
\begin{equation} \label{log}
\phi_{k}=0 \quad \mbox{mod}\,  2\pi
\end{equation}

\section{Expression for the Norm in the XXX Type Models}
For models with an R matrix (\ref{xxxel}) the scalar product of
an eigenfunction  and dual eigenfunction  is equal to
\begin{eqnarray} {\label{xxxnorm}}
& <0|C(\lambda_{1})\ldots C(\lambda_{N}) B(\lambda_{1})\ldots B(\lambda_{N})  |0>= \kappa^{N} \left( \prod_{j=1}^{N}a(\lambda_{j})d(\lambda_{j}) \right) <0|0> \\
& \left( \prod_{j=1}^{N}\prod_{k=1, j\neq k}^{N}  f(\lambda_{j} , \lambda_{k}) \right) \det_{N} \left( \frac{\partial  \phi_{j}} {\partial \lambda_{k}} \right) \nonumber
\end{eqnarray}
Here the $\phi_{k}$ are the variables (\ref{lhs}) and the set of the $\lambda_{j}$  is a solution of the system
(\ref{log}), (\ref{be}).The derivatives can be written in the explicit form:
\begin{equation} \label{dv}
\frac{\partial \phi_{k}}{\partial \lambda_{j}}=\delta_{k,j} \left(i\frac{\partial}{\partial \lambda_{k}}\ln \frac{a(\lambda_{k})}{d(\lambda_{k})}    +\sum_{p=1}^{N} \frac{2\kappa}{(\lambda_{k}-\lambda_{p})^{2}+\kappa^{2}} \right) -\frac{2\kappa}{(\lambda_{k}-\lambda_{j})^{2}+\kappa^{2}}
\end{equation}
The formula (\ref{xxxnorm}) is useful for calculation of norms if the set of $\lambda_{j}$ invariant under complex conjugation
$\{\lambda_{j} \} =\overline{ \{\lambda_{j} \}}$ and $B^{\dagger}(\lambda)=\pm C(\overline{\lambda})$.
Examples of applications can be useful. 
Let consider nonlinear Schroedinger equation, it is also known as Bose gas with delta interaction. It has a Hamiltonian
\begin{eqnarray} \label{contham}
& H=  \int dx  \left( \partial_{x} \psi^{\dagger}  \partial_{x} \psi +\kappa \psi^{\dagger} \psi^{\dagger} \psi \psi \right), \nonumber \\
&  \left[ \psi (x), \psi^{\dagger} (y)  = \delta (x-y)\right]
\end{eqnarray} 
 The monodromy matrix has the following property:
 \begin{equation} \label{nsinv}
   \sigma^{1}T^{\dagger}(\overline{\lambda}) \sigma^{1}= T(\lambda),   \qquad B^{\dagger}(\lambda)=C(\overline{\lambda})
 \end{equation}
  This is $XXX$ case. The vacuum eigenvalues are equal to:
  $$a(\lambda)=\exp (-iL\lambda/2), \qquad  d(\lambda)=\exp (iL\lambda/2)$$
  This proves Gaudin conjecture \cite{gaudin}. This determinant formula for norm also applicable for Heisenberg XXX spin chain.
  
  \section{Expression for the Norm in the $XXZ$ Type Models}
 For solvable models with XXZ $R$ matrix (\ref{xxzel}) the scalar
product of the eigenfunction  and dual eigenfunction  is equal 
  \begin{eqnarray} {\label{xxznorm}}
& <0|C(\lambda_{1})\ldots C(\lambda_{N}) B(\lambda_{1})\ldots B(\lambda_{N})  |0>= (\sin2\eta )^{N} \left( \prod_{j=1}^{N}a(\lambda_{j})d(\lambda_{j}) \right) <0|0> \\
& \left( \prod_{j=1}^{N}\prod_{k=1, j\neq k}^{N}  f(\lambda_{j} , \lambda_{k}) \right) \det_{N} \left( \frac{\partial  \phi_{j}} {\partial \lambda_{k}} \right) \nonumber
\end{eqnarray}
Here $\lambda_{j}$ has to satisfy equation	(\ref{be}). The Jacobi  matrix can be written down in the explicit form:

\begin{equation} \label{dvxxz}  \nonumber
\frac{\partial \phi_{k}}{\partial \lambda_{j}}=\delta_{k,j} \left(i\frac{\partial \ln {a(\lambda_{k})}/{d(\lambda_{k})} }{\partial \lambda_{k}}  +\sum_{p=1}^{N} \frac{\sin (4\eta)}{\sinh (\lambda_{k}-\lambda_{p}+2i\eta ) \sinh (\lambda_{k}-\lambda_{p}-2i\eta )} \right) +\frac{-\sin (4\eta)}{\sinh (\lambda_{k}-\lambda_{j}+2i\eta ) \sinh (\lambda_{k}-\lambda_{j}-2i\eta )}  \nonumber
\end{equation}
The whole paper \cite{origin} is devoted to the derivation of this formula. Let us calculate the norms for X XZ model. The Hamiltonian of the model is
\begin{equation} \label{xxzham} 
H=\sum_{k=1}^{M} \sigma^{1}_{k}\sigma^{1}_{k+1}+ \sigma^{2}_{k}\sigma^{2}_{k+1}+\cos(2\eta) (\sigma^{3}_{k}\sigma^{3}_{k+1}-1 )
\end{equation}
The model was imbedded into QISM in \cite{fadtakh}. This monodromy matrix has the following property at real $\eta$:
  \begin{equation} \label{xxzinv}
   \sigma^{2}T^{\dagger}(\overline{\lambda}) \sigma^{2}= T(\lambda),   \qquad B^{\dagger}(\lambda)=-C(\overline{\lambda})
 \end{equation}
Pseudo-vacuum is the ferromagnetic state $|0>= \prod_{j=1}^{M} \uparrow_{j}=<0|$. The vacuum eigenvalues are:
$$a(\lambda)=\sinh^{M}(\lambda -i \eta), \qquad  d(\lambda)= \sinh^{M}(\lambda +i \eta)$$

In order to write down the square of the norm it is convenient to introduce:
$$\chi (\lambda , \eta) = \frac{\sin (2\eta)}{\sinh (\lambda -{i \eta}) \sinh (\lambda +{i \eta}  )}$$ 
The formula for the square of the norm is:
 \begin{eqnarray} {\label{xxzNorm}} 
&<0|B^{\dagger}(\lambda_{N}) \ldots B^{\dagger}(\lambda_{1})  B(\lambda_{1}) \ldots B(\lambda_{N})|0> = \sin^{N}(2\eta)
 \left(\prod_{j=1}^{N}\sinh^{M} (\lambda_{j}-i\eta) \sinh^{M} (\lambda_{j}+i\eta) \right) \\
&  \left(\prod_{k>j=1}^{N}  \frac{\sinh(\lambda_j -\lambda_k-2i\eta ) \sinh(\lambda_j -\lambda_k+2i\eta )}{\sinh^{2}(\lambda_j -\lambda_k)} \right) \det_{N}\left[\delta_{kj} \left( M\chi (\lambda_{k}, \eta) -\sum_{l=1}^{N}\chi (\lambda_{k}- \lambda_{l}, 2\eta)\right)+ \chi (\lambda_{k}- \lambda_{j}, 2\eta)    \right] \nonumber
\end{eqnarray}
The formula was presented in \cite{gmt} and verifies for $N=2$ and $N=3$. The proof for arbitrary $N$ was first published in \cite{origin}. The formula (\ref{xxznorm}) also describes norms in Sine-Gordon and lattice Sine-Gordon \cite{izko2, izko3}.

\section{The Idea of the Proof}

In order to prove determinant formula for the norm of Bethe wave function let us introduce an object:
\begin{eqnarray} {\label{object}}
 \frac{<0|C(\lambda_{1})\ldots C(\lambda_{N}) B(\lambda_{1})\ldots B(\lambda_{N})  |0>}{(\sin2\eta )^{N}\left( \prod_{j=1}^{N}a(\lambda_{j})d(\lambda_{j}) \right) <0|0>\left( \prod_{j=1}^{N}\prod_{k=1, j\neq k}^{N}  f(\lambda_{j} , \lambda_{k}) \right)}=  |\lambda_{1}\ldots \lambda_{N} | 
\end{eqnarray}
We have to prove that 
\begin{equation} \label{reform}
  |\lambda_{1}\ldots \lambda_{N}|=
 \det_{N} \left( \frac{\partial  \phi_{j}} {\partial \lambda_{k}} \right) 
\end{equation}
We will assume that $<0|0>=1$. The author of \cite{origin} proved that $a(\lambda)$ and $d(\lambda)$ can be considered an arbitrary functions. So we can consider the variables 
$$ X_{p}= i\frac{\partial}{\partial \lambda_{p}} \ln \frac{a(\lambda_{p})}{d(\lambda_{p})}  $$
as independent of $\lambda_{j}$, see \cite{origin} .
The following theorem is proved in  \cite{origin}.

{\it In order to prove (\ref{reform}) in is enough to prove the following five properties of $ |\lambda_{1}\ldots \lambda_{N}|$:
\begin{itemize}
\item It is invariant under simultaneous replacement
\begin{equation} \label{sym}
\lambda_{j} \leftrightarrow \lambda_{k}, \qquad X_{j}  \leftrightarrow X_{k}
\end{equation}
\item
It is leaner function of $X_{1}$
\begin{equation} \label{linear}
|\lambda_{1}\ldots \lambda_{N}| = U_{1}X_{1} +V_{1}
\end{equation}
\item The coefficient $U_{1}$ is 
\begin{equation} \label{modification}
U_{1}= |\lambda_{2}\ldots \lambda_{N}| ^{\mbox{modified}}
\end{equation}
The right hand side is given by formula (\ref{object}) with $\lambda_{1}$ removed and functions $a(\lambda)$ and $d(\lambda)$ replaced by:
\begin{equation} \label{vacmod}
a^{\mbox{modified}}(\lambda)=a(\lambda)f(\lambda , \lambda_{1}), \qquad d^{\mbox{modified}}(\lambda)=d(\lambda)f(\lambda_{1} , \lambda)
\end{equation}
\item It vanish if all $X_{p}=0$
\begin{equation} \label{vanish}
|\lambda_{1}\ldots \lambda_{N}| = 0 , \qquad \mbox{if} \quad \mbox{all} \qquad X_{p}=0, \qquad  \mbox{at} \qquad  p=1, \ldots , N
\end{equation}
\item For $N=1$
$$|\lambda_{1}|=X_{1}$$
\end{itemize} }
The paper \cite{origin} proves that the right hand side of (\ref{object}) has all five properties listed above. The proof is reduced to analysis of six vertex model with domain wall boundary conditions.

\section{Six Vertex Model}
The six vertex model is an important 'counterexample' of statistical mechanics: the bulk free energy depends on the boundary conditions even in thermodynamic limit, see  \cite{kari}.  

Let us start the formal presention:
$L$ operator in site number $k$ is
 \[L_{k}(\lambda- \nu_{k})= \left( \begin{array}{cc}
\sinh(\lambda -\nu_{k} -i\eta \sigma_{k}^{3}) & -i\sigma^{-}_{k}\sin (2\eta)  \\
-i\sigma^{+}_{k}\sin (2\eta)& \sinh(\lambda -\nu_{k} +i\eta \sigma_{k}^{3})  \end{array} \right)\] 
It also can written as:
\begin{equation} \label{Lxxz}
 L_{k}(\lambda- \nu_{k})=\cos(\eta) \sinh(\lambda -\nu_{k}) -i \sin(\eta)\sigma^{3}\sigma^{3}_{k} \cosh(\lambda-\nu_{k}) -i \sin(2\eta) (\sigma^{+}\sigma^{-}_{k}+ \sigma^{-}\sigma^{+}_{k} )
\end{equation}
The monodromy matrix:
\begin{equation} \label{Monxxz}
T(\lambda)=L_{M}(\lambda -\nu_{M}) \dots L_{1}(\lambda -\nu_{1}) 
\end{equation}
obeys the commutation relations (\ref{rtt}) with $R$ matrix (\ref{xxzel}). Pseudo-vacuum is the ferromagnetic state
 $$|0>= \prod_{j=1}^{M} \uparrow_{j}=<0|$$ The vacuum eigenvalues are now equal to
$$a(\lambda)=\prod_{j=1}^{M} \sinh(\lambda -\nu_{j} -i \eta), \qquad  d(\lambda)= \prod_{j=1}^{M} \sinh(\lambda -\nu_{j} +i \eta)$$
Such an inhomogeneous generalization was used for example in \cite{baxter2,bel,kul}. Let us consider a special state
\begin{equation} \label{flip}
B(\lambda_{1})\ldots B(\lambda_{M}) |0>
\end{equation}
Here the number of the $B(\lambda)$ is equal to the number of the sites in
the lattice $N=M$. All spins are looking down in this state:
\begin{equation} \label{afer}
B(\lambda_{1})\ldots B(\lambda_{M}) |0>= Z_{N} \Omega
\end{equation}
Here $Z_{N}$ is a complex number and
\begin{equation} \label{omega}
\Omega= \prod_{j=1}^{M} \downarrow_{j}
\end{equation}
The definition of $Z_{N}$ is
$$Z_{N}= \Omega B(\lambda_{1})\ldots B(\lambda_{M}) |0>$$
The paper \cite{origin} proves that $Z_{N}$ is the partition function of six vertex model with domain wall boundary conditions, see Appendix D.
Actually these boundary conditions were introduced in this paper, including the name. Explicit description of six vertex model with domain wall boundary conditions also can be found in \cite{zinn}.

The following recursion relations were discovered in \cite{origin}:
If $\nu_{1}=\lambda_{1}+i\eta$ then $Z_{N}$ reduces to $Z_{N-1}$ with $\lambda_{1}$ and $\nu_{1}$ removed.
 \begin{eqnarray} {\label{rec}} \nonumber
 Z_{N}(\{\lambda_{\alpha} \}, \{ \nu_{j} \} )|_{\nu_{1}-\lambda_{1}=i\eta}=-i\sin(2\eta) \left(  \prod_{k=2}^{N} \sinh(\lambda_{1} -\nu_{k} -i \eta) \right)  \left(  \prod_{\alpha=2}^{N} \sinh(\lambda_{\alpha} -\nu_{1} -i \eta) \right)  Z_{N-1}(\{\lambda_{\alpha \neq 1} \}, \{ \nu_{j \neq 1} \} )
 \end{eqnarray} 
 The  derivation of this recursion relations also  can  be found in section 2 of  the paper \cite{zinn}.
 
  \section{Conclusion}
 There was a lot of progress since the paper \cite{origin}.  The partition function of six vertex model with domain wall boundary conditions has a lot of applications \footnote{Goolge finds thousands of publications for: Six vertex model with domain wall boundary conditions } and generalizations, see for example \cite{rose}.  Still there are open problems. For example the determinant formula for Bethe wave function in Hubbard still not proven, see \cite{hubbard}.


\begin{thebibliography}{99}
\bibitem{beth} Bethe,H .: Z. Phys.71, 205-226( 1931)  
\bibitem{yy} Yang.C .N., Yang,C .P.: Phys.Rev.150, 321-327( 1966)
\bibitem{lieb} Lieb.E .H. : Phys. Rev. Lett . 18. 692-694 ( 1967)
\bibitem{bpf} Berezin F. A,  .Pokhil ,C .P,  Finkelberg V, M. : Vestnk Mosk. Gos.Univ.Ser 1.1, 21-28( 1964) 
\bibitem{mcg} McGuire J. , B.   J.Math. Phys 5,  622-636(1964)

Brezin.E , Zinn-Justin J, : C.R. Acad. Sci.Paris,  263, 670- 613 (1966)

Gaudin, M ,  : J. Math. Phys. 12, 1674-1680, (1971)
\bibitem{fadd} Faddeev L, D. :  Sov. Sci. R ev. Math. Phys C. l , 107-160 (1981)
\bibitem{fadtakh} Faddeev L,. D. , Takhtajan L, . A. : Usp. Mat . Nauk 34, 13- 63 (1979)
\bibitem{gaudin} Gaudin. M . :  Preprint Centre d 'Etudes Nucleaire de Saclay CEA-N-1559 (1) , ( 1972)
\bibitem{gmt} Gaudin, M . , McCoy, B .M. , Wu,T .T.  :Phys. Rev.D, vol  23, 417 (1981)
\bibitem{baxter}  Baxter R, . J. : J. Stat . Phys. 9, 145 -182 ( 1973)
\bibitem{kari}  Eloranta K., J. Stat. Phys. 96 (1999), 1091.
\bibitem{ik} Izergin A, .G. , Korepin V, .E. : L ett . Math. Phys. (to be published)
\bibitem{fasklya} Faddeev L,. D. , Sklyanin E, .K. : Dokl . Akad. N auk SSSR 243, 1430-1433 (1978)
\bibitem{slya} Sklyanin E, .K. : Dokl . Akad. N auk SSSR 244, 1337-1341 (1978)
\bibitem{iks} Izergin A, . G. , Korepin V, .E. , Smirnov F, .A. : Teor . Mat . Fiz. 48, 319-323 (1981)
\bibitem{slkya1} Sklyanin E, .K. : Zap. Nauchn Seminarov LOMI,  95, 55-128 (1980)
\bibitem{izko1} Izergin A, .G. , Korepin V, .E. : Dokl . Akad. N auk SSSR 259, 76-79 (1981)
\bibitem{izko2}  Izergin. A.G.,  Korepin, V.E. :  Nucl. Phys. B. Field Theory and Statistical Systems B205 [FS5]
401-413 (1982)
\bibitem{izko3} Izergin, A.G..  Korepin, V.E. : Lett. Math. Phys. 5, 199 -205 (1981)
\bibitem{baxter2} Baxter, R.J.: Stud. Appl. Math. L50, 51-67 (1971)
\bibitem{bel} Belavin. A.A.:  Phys. Lett. B87, 117-121 (1980)
\bibitem{kul} Kulish. P.P. : Physica  D3, 246-257 (1981)
\bibitem{fadta3} Faddeev, L.D.,  Takhtajan, L.A.: Zap. Nauchn. Seminarov LOMI 109, 134-178 (1981)
\bibitem{baxter5} Baxter, R.J.: Philos. Trans. Soc. London, Ser. A  289, 315-346 (1978)
\bibitem{origin} Korepin V.E,  Communications in Mathematical Physics, vol 86, page 391 in l982.

To see the paper  you can go to the  page http://insti.physics.sunysb.edu/physics/forms/profilesearch.cgi?lastname=korepin
 click PERSONAL HOMEPAGE (in the right column) and find the .PDF file for Calculation of Norms of Bethe Wave Functions (fifth bullet from above).
\bibitem{book} V.E. Korepin, N.M. Bogoliubov and A.G. Izergin, a book Quantum Inverse Scattering Method and Correlation Functions, Cambridge University Press, 1993 , http://www.cambridge.org/catalogue/catalogue.asp?isbn=9780521586467
\bibitem{hubbard} F. Goehmann, V. E. Korepin, Phys.Lett. A263 (1999) 293-298, arXiv:cond-mat/9908114 
\bibitem{zinn} V. Korepin, P. Zinn-Justin, arXiv:cond-mat/0004250 , J. Phys. A 33 No. 40 (2000), 7053
\bibitem{rose} Hjalmar Rosengren http://arxiv.org/abs/0911.0561 
\end{thebibliography}
\end{document}